# Data-driven Topology Optimization (DDTO) for Three-dimensional Continuum Structures

Yunhang Guo[1], Zongliang Du[1,2]*, Lubin Wang[2], Wen Meng[3], Tien Zhang[3], Ruiyi Su[3], Dongsheng Yang[4], Shan Tang[1,2]*, Xu Guo[1,2]*

[1]*State Key Laboratory of Structural Analysis for Industrial Equipment,*
*Department of Engineering Mechanics,*
*Dalian University of Technology, Dalian, 116023, P.R. China*
[2]*Ningbo Institute of Dalian University of Technology, Ningbo, 315016, P.R. China*
[3]*Beijing System Design Institute of Electro-Mechanic Engineering, Beijing, 100854, P.R. China*
[4]*China Academy of Launch Vehicle Technology, Beijing, 100076, P.R. China*

**Abstract**

Developing appropriate analytic-function-based constitutive models for new materials with nonlinear mechanical behavior is demanding. For such kinds of materials, it is more challenging to realize the integrated design from the collection of the material experiment under the classical topology optimization framework based on constitutive models. The present work proposes a mechanistic-based data-driven topology optimization (DDTO) framework for three-dimensional continuum structures under finite deformation. In the DDTO framework, with the help of neural networks and explicit topology optimization method, the optimal design of the three-dimensional continuum structures under finite deformation is implemented only using the uniaxial and equi-biaxial experimental data. Numerical examples illustrate the effectiveness of the data-driven topology optimization approach, which paves the way for the optimal design of continuum structures composed of novel materials without available constitutive relations.

*Corresponding authors. E-mail: zldu@dlut.edu.cn (Zongliang Du), shantang@dlut.edu.cn (Shan Tang), guoxu@dlut.edu.cn (Xu Guo)

# 1. Introduction

Structural topology optimization aims to find the optimal material distribution under prescribed objective functions and constraint conditions in the design domain. Since the pioneering work of Bendsøe and Kikuchi [1], topology optimization of the continuum structure has received a lot of interest from researchers and engineers. Nowadays, many structural topology optimization methods have been proposed, such as the solid isotropic material with penalization (SIMP) method [2, 3], the level set method (LSM) [4, 5], the evolutionary structural optimization (ESO) method [6], the moving morphable components (MMC) method [7], etc. Although significant progress has been achieved in the optimal design of continuum structures, most research is carried out under the linear elasticity and small deformation assumptions, which inevitably restrict the applications.

In terms of the topology optimization of continuum structures with material and geometric nonlinearities, some excellent progress has been reported recently. Based on the LSM, Chen et al. [8] designed hyperelastic structures undergoing large deformation. Dalklint et al. [9] studied the eigenfrequency constrained topology optimization of hyperelastic structures. Deng et al. [10] presented a topology optimization method based on distortion energy for designing hyperelastic material against failure. Recently, topological design of porous infill structures with hyperelastic material under large deformation configurations was reported in [11]. Using a modified evolutionary topology optimization method, Zhang et al. [12] maximized the stiffness of hyperelastic structures in the finite deformation regime. Silva et al. [13] proposed a modified Normal Distribution Fiber Optimization (NDFO) model to perform topology optimization of fibers orientation in composites. Based on the hyperelasticity theory, Zhang et al. [14, 15] implemented the structural topology optimization design of multi-material and local feature control, respectively. It is worth noting that various techniques which alleviate the convergence difficulties caused by the excessive distortion of low-density elements, play a vital role in the topology optimization of continuum structure with geometric nonlinearity, e.g., see [16-25] for reference.

It is worth mentioning that, all of the studies mentioned above were investigated under the constitutive model-based structural topology optimization framework. On the one hand, developing a matured analytic function-based constitutive model, in general, is very demanding; on the other hand, there may be no available accurate constitutive models for new materials, such as many synthetic polymers and additive manufactured materials. To be specific, reference [26] shows that the neo-Hookean model, Arruda–Boyce model, and Odgen model cannot capture the mechanical property accurately under compression and/or tension for fabricated silicone rubber, see Fig. 1. At this circumstance, it would be extremely difficult to design optimal structures only with the collection of experimental data of such materials. Fortunately, Kirchdoerfer and Ortiz [27] proposed a data-driven computational mechanics (DDCM) framework, which realized structural mechanics analysis based only on experimental material data. Later works suggested performing the DDCM based on locally characterizing the material constitutive behavior to accelerate the global data-driven solver and improve the robustness against noise [28-31]. To consider the inevitable multi-source uncertainties in the material data set, Guo et al. [32] proposed an uncertainty analysis-based data-driven computational mechanics (UA-DDCM) framework to obtain a solution set rather than a specific value of the concerned structural response. Tang et al. [26, 33-35] developed a so-called MAP123 method based on the principle of mechanics. They successfully realized the mechanical analysis of the three-dimensional structures only using experimental material data of uniaxial tension, or uniaxial and equi-biaxial tension.

Inspired by the works in data-driven computational mechanics, the idea of data-driven topology optimization (DDTO) emerges naturally [36], which aims at designing optimal structures only using experimental material data, or in other words, based on the constitutive-model-free data-driven structural analysis method. To the best of our knowledge, however, only one related work was reported by Zhou et al. [36] on this topic. In that work, combined with the one-dimensional experimental material data, topology optimization of the truss structures is performed as a three-level programming based on the DDCM method under the assumption of small deformation. Besides the novel and encouraging results achieved in that work, there are still some challenging issues when the DDTO is thought for three-dimensional continuum structures under finite deformation: (1) From the data-driven analysis aspect, due to the existence of complex stress states in the intermediate three-dimensional continuum designs under finite deformation, a stable and efficient data-driven structural analysis method is required. In addition, how to effectively construct and utilize the data set for a three-dimensional continuum structure is another difficulty in the DDCM framework. (2) From the topology optimization aspect, even efficient data-driven analysis algorithm is available, the mathematical formulation, sensitivity analysis results, and numerical stabilization techniques (e.g., alleviation of the excessive distortion of weak elements) in the DDTO framework desire further exploration.

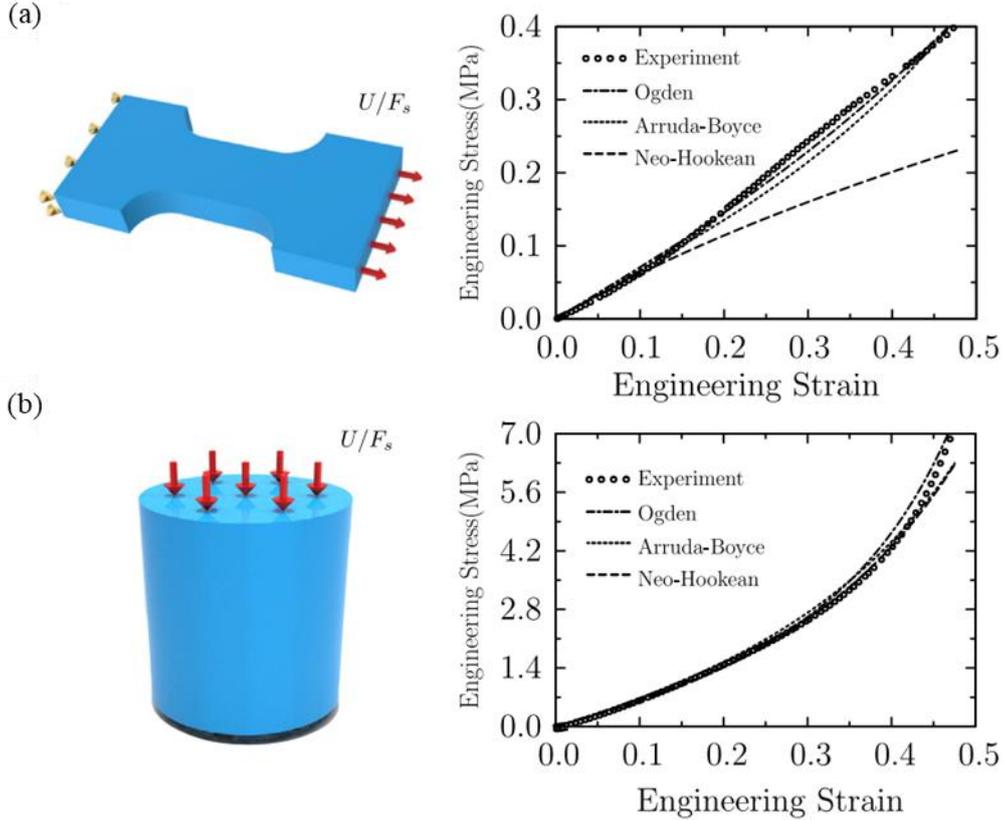

Fig. 1. Stress–strain data for fabricated silicone rubber in experiments under (a) tension and (b) compression respectively. All of the neo-Hookean model, Arruda–Boyce model or Odgen model cannot perfectly fit the experimental data when the deformation is large [26].

In the present work, a mechanistic-based data-driven topology optimization (DDTO) framework for three-dimensional continuum structures is proposed to address the above-mentioned challenging issues. Under the assumptions of isotropy and the coaxial relationship between deviatoric stress and deviatoric strain, the MAP123 method is improved with the help of artificial neural networks, which can predict a hyperelastic constitutive model locally by only utilizing the uniaxial and equi-biaxial experimental data. Based on this, the DDTO framework can be effectively tackled with the help of the theory and solution techniques of the topology optimization of three-dimensional hyperelastic structures.

The remainder of the paper is organized as follows. In Section 2, the classical constitutive model-based topology optimization framework for three-dimensional hyperelastic continuum structures under finite deformation and the corresponding solution strategies are explained in detail. In Section 3, after developing a stable data-driven analysis method for three-dimensional continuum structures, the DDTO framework, its sensitivity analysis and solution strategy are presented. After numerical examples illustrating the effectiveness of the DDTO framework in Section 4, some concluding remarks are provided in Section 5.

## 2. Constitutive model-based explicit topology optimization framework for three-dimensional hyperelastic continuum structures under finite deformation

For completeness, this section mainly introduces the classical constitutive model-based topology optimization (CMTO) framework for three-dimensional hyperelastic continuum structures, including subsections about the three-dimensional Moving Morphable Void (MMV) method, problem formulation, and numerical solution strategies.

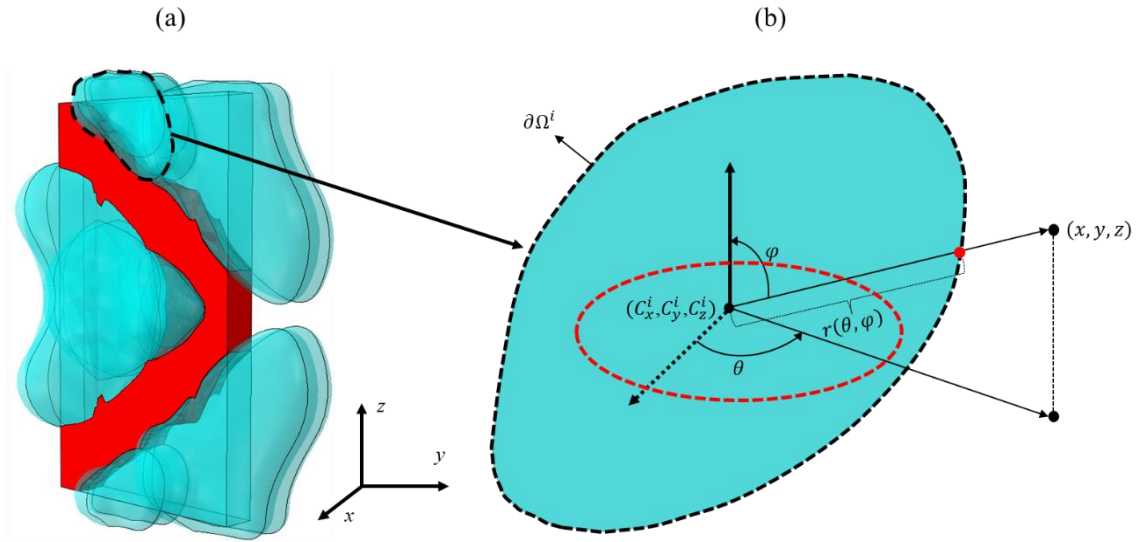

Fig. 2. The schematic diagram of the 3D-MMV method, (a) an optimized structure, (b) geometry representation of the $i$th void.

### 2.1 Three-dimensional Moving Morphable Void (MMV) method

To solve the topology optimization problems of linear elastic three-dimensional continuum structures, Zhang et al. [37] proposed the so-called 3D-MMV method. As shown in Fig. 2(a), an optimized structure is described by a set of movable and morphable voids $\Omega_1, \ldots, \Omega_{nv}$ constructed by NURBS surfaces, where $nv$ is the total number of voids in the design domain. The $i$th void is determined by its topology description function (TDF) $\chi^i(\boldsymbol{x})$ defined as

$$\chi^i = \sqrt{\left(x - C_x^i\right)^2 + \left(y - C_y^i\right)^2 + (z - C_z^i)^2} - r(\theta, \varphi), \qquad (2.1)$$

where $\left(C_x^i, C_y^i, C_z^i\right)^\top$ is the coordinate of its center, as illustrated by Fig. 2(b), $\theta$ and $\varphi$ are the azimuths, $r(\theta, \varphi)$ is the distance function from the point on the boundary $\partial \Omega^i$ to the center point of the void and its detailed expression is listed in Appendix A.

Based on Eq. (2.1), the global TDF of the whole structure can be constructed as [38]:

$$\chi^s(\boldsymbol{D}; \boldsymbol{x}) = \min_{i=1,\ldots,nv} \chi^i(\boldsymbol{D}^i; \boldsymbol{x}) \approx \frac{1}{-\lambda} \ln \left( \sum_{i=1}^{nv} e^{-\lambda \chi^i(\boldsymbol{D}^i; \boldsymbol{x})} \right), \qquad (2.2)$$

where $\boldsymbol{D}^i = \left(C_x^i, C_y^i, C_z^i, r^{0,0}, r^{0,L}, \ldots, r^{k,l}, \ldots, r^{K,L-1}\right)^\top$, $k = 1, \ldots, K$, $l = 1, \ldots, L-1$, and $\lambda$ is a relatively large even number, e.g., $\lambda = 50$. Notably, the design variable vector in the MMV approach is $\boldsymbol{D} = ((\boldsymbol{D}^1)^\top, \ldots, (\boldsymbol{D}^{nv})^\top)^\top$ and it only contains geometrically explicit parameters.

**2.2 Problem formulation**

By adopting the total Lagrangian formulation, for minimum end compliance design problems, the constitutive model based-explicit topology optimization framework for three-dimensional hyperelastic continuum structures can be formulated as:

Find $\quad \boldsymbol{D} = \left((\boldsymbol{D}^1)^\top, \ldots, (\boldsymbol{D}^i)^\top, \ldots, (\boldsymbol{D}^{nv})^\top\right)^\top$, $^0\boldsymbol{u}(\boldsymbol{D}; {}^0\boldsymbol{x})$,

Minimize $\quad I(\boldsymbol{D}, {}^0\boldsymbol{u}) = \int_{{}_0\Omega} H_\epsilon^\alpha\left(\chi^s(\boldsymbol{D}; {}^0\boldsymbol{x})\right) {}^0\boldsymbol{f} \cdot {}^0\boldsymbol{u} \mathrm{d}^0\mathrm{V} + \int_{{}_0S_t} {}^0\boldsymbol{t} \cdot {}^0\boldsymbol{u} \mathrm{d}^0\mathrm{S}$,

S. t.

$$\int_{{}_0\Omega} H_\epsilon^\alpha\left(\chi^s(\boldsymbol{D}; {}^0\boldsymbol{x})\right) \delta \boldsymbol{F}^\top : \boldsymbol{P} \mathrm{d}^0\mathrm{V} = \int_{{}_0\Omega} H_\epsilon^\alpha\left(\chi^s(\boldsymbol{D}; {}^0\boldsymbol{x})\right) {}^0\boldsymbol{f} \cdot \delta {}^0\boldsymbol{u} \mathrm{d}^0\mathrm{V} + \int_{{}_0S_t} {}^0\boldsymbol{t} \cdot \delta {}^0\boldsymbol{u} \mathrm{d}^0\mathrm{S}, \forall \delta {}^0\boldsymbol{u} \in \mathcal{U}_{\mathrm{ad}}^0,$$

$$\int_{^0\Omega} H_\epsilon^\alpha\bigl(\chi^s(\boldsymbol{D};\,^0\boldsymbol{x})\bigr)\mathrm{d}^0\mathrm{V} \leq {^0\bar{V}},$$

$$\boldsymbol{D} \subset \mathcal{U}_D,$$

$$^0\boldsymbol{u} = \boldsymbol{0},\,\mathrm{on}\,^0\mathrm{S_u}, \qquad (2.3)$$

where $I(\boldsymbol{D},\,^0\boldsymbol{u})$ is the end compliance of the optimized structure, $^0\Omega$ is the undeformed design domain and $^0(\ )$ denotes the variables in the initial configuration. The 1st Piola–Kirchhoff stress tensor $\boldsymbol{P}$ is obtained from the Cauchy stress tensor $\boldsymbol{\sigma}$ determined by the strain energy density function $W$ (see Section 3), i.e., $\boldsymbol{P} = J\boldsymbol{\sigma}\cdot\boldsymbol{F}^{-\mathrm{T}}$ with $\boldsymbol{F}$ and $J$ denoting the deformation gradient and the determinant of $\boldsymbol{F}$. Besides, $^0\boldsymbol{f}$ is the body force density, $^0\boldsymbol{u}$ is the displacement field of the structure, $\delta^0\boldsymbol{u}$ is virtual displacement, $^0\boldsymbol{t}$ is surface tractions density, $\mathcal{U}_{\mathrm{ad}}^0$ is the admissible set that $\delta^0\boldsymbol{u}$ belongs to, $^0\mathrm{S_u}$ is the Dirichlet boundary, $^0\bar{V}$ is the upper bound of admissible solid material. The regularized Heaviside function $H_\epsilon^\alpha$ is expressed as

$$H_\epsilon^\alpha(x) = \begin{cases} 1, & \text{if } x > \epsilon, \\ \dfrac{3(1-\alpha)}{4}\left(\dfrac{x}{\epsilon} - \dfrac{x^3}{3\epsilon^3}\right) + \dfrac{1+\alpha}{2}, & \text{if } -\epsilon \leq x \leq \epsilon, \\ \alpha, & \text{otherwise}, \end{cases} \qquad (2.4)$$

where $\epsilon$ is the width of regularized region and $\alpha$ is a small positive number, which are set as $0.1$ and $10^{-3}$, respectively. In the solution process of formulation (2.3), the strain energy density function $W_e$ of the hyperelastic material corresponding to the $e$th element is obtained according to the ersatz material model as

$$W_e = \rho_e W_0 = \frac{\sum_{i=1}^{nn} H_\epsilon^\alpha\bigl(\chi^s(\boldsymbol{D};\,^0\boldsymbol{x})_i^e\bigr)}{nn} W_0, \qquad (2.5)$$

where $\rho_e$ is the density of the $e$th element, $nn$ is the total number of the nodes in $e$th element, $\chi^s(\boldsymbol{D};\,^0\boldsymbol{x})_i^e$ is the $i$th nodal value of the global TDF of the $e$th element, $W_0$ is the strain energy density functional of the base material.

**2.3 Numerical solution strategies**

*2.3.1 Redundant degrees of freedom removal technique*

The solution process of topology optimization of three-dimensional hyperelastic continuum structures under finite deformation usually encounters the non-convergence issue of the finite element analysis caused by the excessive deformation of weak material elements. The so-called redundant degrees of freedom removal technique [25, 37, 38] is adopted to alleviate this issue. The core idea of this technique is to remove the degrees of freedom (DOFs) only involved in the weak material elements which are identified by Eq. (2.6) for the nonlinear structural analysis:

$$\rho_e = \frac{\sum_{i=1}^{nn} H_\epsilon^\alpha(\chi^s(\bm{D}; {}^0\bm{x})_i^e)}{nn} \leq TR, \qquad (2.6)$$

where $TR$ is a threshold taken as 0.1 in the present work[1]. It is worth noting that, thanks to the decoupling between the geometry description and finite element analysis, one would not need to worry about the re-introduction of the removed DOFs in subsequent iterations. For some more detailed explanations, please refer to [25, 37, 38].

*2.3.2 Finite element analysis*

Combining the weak form of the balance equation in Eq. (2.3) and the redundant degrees of freedom removal technique, the balance equation of the retained degrees of freedom in the design domain can be written as:

$$\widetilde{\bm{R}}(\bm{D}, \widetilde{\bm{U}}) = \widetilde{\bm{f}}^{\text{int}} - \widetilde{\bm{f}}^{\text{ext}} = \sum_{e=1}^{Rne} \bm{L}_e^\top \bm{f}_e^{\text{int}} - \sum_{e=1}^{Rne} \bm{L}_e^\top \bm{f}_e^{\text{ext}} = \bm{0}, \qquad (2.7)$$

where $\widetilde{\bm{U}}$, $\widetilde{\bm{R}}(\bm{D}, \widetilde{\bm{U}})$, $\widetilde{\bm{f}}^{\text{int}}$ and $\widetilde{\bm{f}}^{\text{ext}}$ are the displacement vector, the residual force vector, the internal force vector and the external force vector of the retained DOFs in the design domain, respectively. The symbol $Rne$ is the number of the elements with density larger than the threshold in the design domain, and $\bm{L}_e$ is the connection matrix of the $e$th element. The internal force vector $\bm{f}_e^{\text{int}}$ and the external force vector $\bm{f}_e^{\text{ext}}$ of the $e$th element can be expressed as:

---

[1]It will be shown by numerical examples, this threshold could balance the numerical stability and accuracy. And we note that, for topology optimization of 2D nonlinear structures, a smaller threshold (e.g., 0.01) can be chosen [25]. This implies the convergence issue caused by excessive deformation of weak material elements in 3D case is more severe.

$$\boldsymbol{f}_e^{\text{int}} = \int_{{}^0\Omega_e} \boldsymbol{B}^\top \boldsymbol{P}\, \mathrm{d}{}^0V_e, \quad \boldsymbol{f}_e^{\text{ext}} = \int_{{}^0\Omega_e} \boldsymbol{N}^\top {}^0\boldsymbol{f}\, \mathrm{d}{}^0V_e + \int_{{}^0S_t^e} \boldsymbol{N}^\top {}^0\boldsymbol{t}\, \mathrm{d}{}^0S, \tag{2.8}$$

where $\boldsymbol{B}$ is the operation matrix between $\delta \boldsymbol{F}_e$ and the virtual displacement $\delta \boldsymbol{U}_e$, $\boldsymbol{N}$ is the shape function matrix. Implementing the Newton-Raphson algorithm, the displacement vector $\widetilde{\boldsymbol{U}}^{k+1}$ at the $k$th iteration is updated as

$$\begin{cases} \boldsymbol{K}_{\text{T}}^k(\boldsymbol{D}, \widetilde{\boldsymbol{U}}^k) \Delta \widetilde{\boldsymbol{U}}^k = -\widetilde{\boldsymbol{R}}(\boldsymbol{D}, \widetilde{\boldsymbol{U}}^k) = \widetilde{\boldsymbol{f}}^{\text{ext}} - \widetilde{\boldsymbol{f}}^{\text{int}}(\boldsymbol{D}, \widetilde{\boldsymbol{U}}^k), \\ \widetilde{\boldsymbol{U}}^{k+1} = \widetilde{\boldsymbol{U}}^k + \Delta \widetilde{\boldsymbol{U}}^k, \end{cases} \tag{2.9}$$

where $\boldsymbol{K}_{\text{T}}^k(\boldsymbol{D}, \widetilde{\boldsymbol{U}}^k)$ and $\Delta \widetilde{\boldsymbol{U}}^k$ are the tangent stiffness matrix and the nodal displacement at the $k$th iteration. If $\|\widetilde{\boldsymbol{R}}(\boldsymbol{D}, \widetilde{\boldsymbol{U}}^k)\|/\|\widetilde{\boldsymbol{f}}^{\text{ext}}\| < \delta$ with $\delta$ denoting a prescribed small number, e.g., $\delta = 10^{-6}$, the structure is supposed to achieve an equilibrium state.

*2.3.3 Sensitivity analysis*

Without the loss of generality, the body force is ignored in the present work. Based on the adjoint sensitivity analysis method, the Lagrangian function only involving the retained DOFs can be written as:

$$L(\boldsymbol{D}, \widetilde{\boldsymbol{U}}, \widetilde{\boldsymbol{\lambda}}) = \widetilde{\boldsymbol{f}}^{\text{ext}\,\top} \widetilde{\boldsymbol{U}} + \widetilde{\boldsymbol{\lambda}}^\top \widetilde{\boldsymbol{R}}(\boldsymbol{D}, \widetilde{\boldsymbol{U}}), \tag{2.10}$$

where $\widetilde{\boldsymbol{\lambda}}$ is the corresponding adjoint displacement vector. With Eq. (2.10) in hand, the derivative of the objective function with respect to the design variables is expressed as:

$$\frac{\partial I(\boldsymbol{D}, {}^0\boldsymbol{u})}{\partial \boldsymbol{D}} = \frac{\partial L(\boldsymbol{D}, \widetilde{\boldsymbol{U}}, \widetilde{\boldsymbol{\lambda}})}{\partial \boldsymbol{D}} = \frac{1}{nn} \widetilde{\boldsymbol{\lambda}}^\top \sum_{e=1}^{Rne} \left( \sum_{i=1}^{nn} \frac{\partial H_\epsilon^\alpha(\chi^s)_i^e}{\partial \chi^s} \frac{\partial \chi^s}{\partial \boldsymbol{D}} \right) \frac{\widetilde{\boldsymbol{f}}_e^{\text{int}}}{\rho_e}. \tag{2.11}$$

The detailed expressions of $\frac{\partial \chi^s}{\partial \boldsymbol{D}}$ are presented in Appendix B and the adjoint displacement vector can be solved by the following adjoint equation:

$$\widetilde{\boldsymbol{\lambda}} = -\left( \boldsymbol{K}_{\text{T}}^E(\boldsymbol{D}, \widetilde{\boldsymbol{U}}^E) \right)^{-1} \widetilde{\boldsymbol{f}}^{\text{ext}}, \tag{2.12}$$

where $\boldsymbol{K}_{\text{T}}^E(\boldsymbol{D}, \widetilde{\boldsymbol{U}}^E)$ and $\widetilde{\boldsymbol{U}}^E$ are the tangent stiffness matrix and the displacement vector at the equilibrium state, respectively. The sensitivity result of the volume constraint is:

$$\frac{\partial V}{\partial \boldsymbol{D}} = \frac{1}{nn} \sum_{e=1}^{ne} \sum_{i=1}^{nn} \frac{V_e}{V_D} \frac{\partial H_\epsilon^\alpha(\chi^s)_i^e}{\partial \chi^s} \frac{\partial \chi^s}{\partial \boldsymbol{D}}, \tag{2.13}$$

where $ne$ is the total number of elements in the structure, the symbols $V_e$ and $V_D$ are the volumes of the $e$th element and the design domain, respectively.

## 3. The DDTO framework for three-dimensional continuum structures under finite deformation

In this section, firstly, a stable data-driven structural analysis (DDSA) algorithm for three-dimensional continuum structures under finite deformation is proposed. Then the effectiveness of DDSA algorithm is verified by a numerical example. Finally, the solution techniques of the corresponding DDTO framework are given.

### 3.1 Neural network enhanced DDSA algorithm for three-dimensional continuum structures under finite deformation

In Tang et al. [26, 35], a constitutive model-free DDSA algorithm for three-dimensional hyperelastic continuum structures is proposed only using the experimental data of uniaxial tension, or uniaxial and equi-biaxial tension tests. In the present work, we also assume the mechanical property of the considered material is approximately hyperelastic. To improve the robustness of the algorithms in [26, 35] for intermediate structures in the optimization process, here, an improved DDSA algorithm, which additionally uses the experimental data of uniaxial compression and equi-biaxial compression test, is proposed. Furthermore, an artificial neural network-based stress update strategy is proposed to ensure the stability and robustness of the structural analysis.

*3.1.1 The prior knowledge in mechanics*

For an isotropic hyperelastic material, its mechanical behavior can be characterized by a strain energy density function expressed as

$$W = W(\bar{I}_1, \bar{I}_2, J), \tag{3.1}$$

where $\bar{I}_1 = \mathrm{tr}(\bar{\boldsymbol{b}})$ and $\bar{I}_2 = \frac{1}{2}\left[\bar{I}_1^{\,2} - \mathrm{tr}(\bar{\boldsymbol{b}} \cdot \bar{\boldsymbol{b}})\right]$ are the first and second invariants of $\bar{\boldsymbol{b}}$, respectively. The symbol $\bar{\boldsymbol{b}} = \bar{\boldsymbol{F}} \cdot \bar{\boldsymbol{F}}^\top$ is the modified left Cauchy-Green tensor and $\bar{\boldsymbol{F}} = J^{-\frac{1}{3}}\boldsymbol{F}$ is the modified deformation gradient with the volume change eliminated. The Cauchy stress can be

derived as:

$$\boldsymbol{\sigma} = \frac{2}{J^{\frac{5}{3}}}\left(\frac{\partial W}{\partial \bar{I}_1} + \bar{I}_1\frac{\partial W}{\partial \bar{I}_2}\right)\boldsymbol{b} - \frac{2}{J^{\frac{7}{3}}}\frac{\partial W}{\partial \bar{I}_2}\boldsymbol{b}^2 + \left(\frac{\partial W}{\partial J} - \frac{2}{3J}\frac{\partial W}{\partial \bar{I}_1}\bar{I}_1 - \frac{4}{3J}\frac{\partial W}{\partial \bar{I}_2}\bar{I}_2\right)\boldsymbol{I}, \qquad (3.2)$$

where $\boldsymbol{I}$ is the second order identity tensor and $\boldsymbol{b} = \boldsymbol{F}\cdot\boldsymbol{F}^\top$ is the left Cauchy-Green tensor. Furthermore, the spherical part and the deviatoric part of the Cauchy stress can be obtained as:

$$\sigma_\mathrm{m} = \frac{\partial W}{\partial J}, \qquad (3.3\mathrm{a})$$

$$\mathrm{dev}(\boldsymbol{\sigma}) = \gamma\,\mathrm{dev}(\boldsymbol{b}) + \beta\,\mathrm{dev}(\boldsymbol{b}^2), \qquad (3.3\mathrm{b})$$

where $\gamma = \frac{2}{J^{\frac{5}{3}}}\left(\frac{\partial W}{\partial \bar{I}_1} + \bar{I}_1\frac{\partial W}{\partial \bar{I}_2}\right)$, $\beta = -\frac{2}{J^{\frac{7}{3}}}\frac{\partial W}{\partial \bar{I}_2}$, and $\mathrm{dev}(\cdot)$ denotes the deviatoric operator.

It is clear that, uniaxial and equi-biaxial experimental data are sufficient for determining the spherical part of the Cauchy stress. However, due to the complexity of the deviatoric part of the Cauchy stress in high dimensions, it is impossible to represent the material property information of the entire state space by only using the uniaxial and equi-biaxial experimental data. This motivates us to reconstruct the deviatoric part of the Cauchy stress with the help of material parameters $\gamma$ and $\beta$.

To determine the material parameters $\gamma$ and $\beta$, the equivalent Cauchy stress $\sigma_\mathrm{e}$ and equivalent strain $b_\mathrm{e}$ are defined as:

$$\sigma_\mathrm{e} = \sqrt{\frac{3}{2}\mathrm{dev}(\boldsymbol{\sigma}):\mathrm{dev}(\boldsymbol{\sigma})}, \qquad (3.4\mathrm{a})$$

$$b_\mathrm{e} = \sqrt{\frac{2}{3}\mathrm{dev}(\boldsymbol{b}):\mathrm{dev}(\boldsymbol{b})}. \qquad (3.4\mathrm{b})$$

Combining Eq. (3.4a), (3.4b) and (3.3b), we have

$$\frac{2}{3}\sigma_\mathrm{e}^2 = \frac{3}{2}\gamma^2 b_\mathrm{e}^2 + 2\gamma\beta\frac{2}{3}\mathrm{dev}(\boldsymbol{b}):\mathrm{dev}(\boldsymbol{b}^2) + \beta^2\mathrm{dev}(\boldsymbol{b}^2):\mathrm{dev}(\boldsymbol{b}^2). \qquad (3.5)$$

Generally, the strain energy density function $W$ of a hyperelastic material includes two cases: 1) $W$ is dependent on $\bar{I}_1$, $\bar{I}_2$ and $J$, e.g., the Mooney-Rivlin and Van der Waals models, and 2) $W$ is only dependent on $\bar{I}_1$, $J$, e.g., the Arruda-Boyce and Yeoh models. For the sake of simplicity, only the latter case is considered in this work, the more complicated case is referred to [35] for reference. Under this circumstance, with $W$ depending on $\bar{I}_1$ and $J$, $\beta$ equals 0, and Eq. (3.5) can be simplified as:

$$\gamma = \frac{2}{3}\frac{\sigma_e}{b_e}. \tag{3.6}$$

Once the material coefficients $\gamma$ is determined, the deviatoric part of the Cauchy stress can be accurately calculated by Eq. (3.3b).

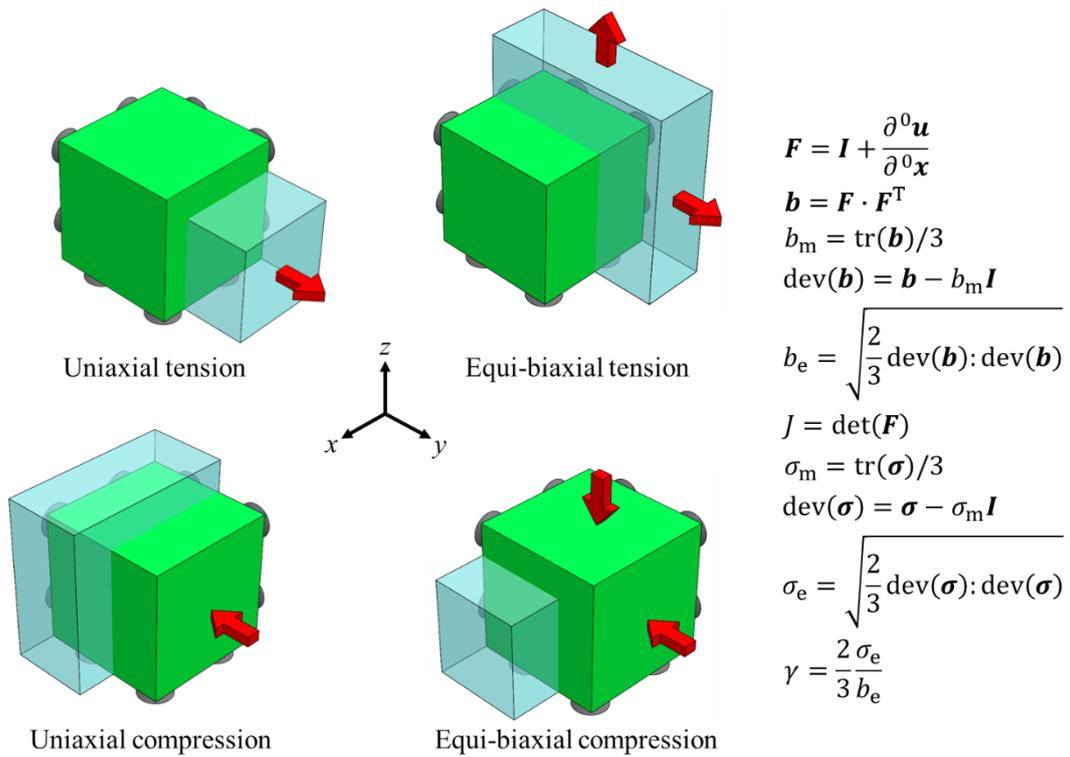

Fig. 3. The schematic diagram of data generation through numerical experiments of uniaxial tension, uniaxial compression, equi-biaxial tension, and equi-biaxial compression.

*3.1.2 Database construction*

In this work, all the experimental data used are generated by numerical experiments of finite element analysis (FEA). Unless otherwise specified, the analytic-function based constitutive model referred to in this work is the Arruda-Boyce model, which is expressed as:

$$W = \mu \sum_{i=1}^{5} \frac{C_i}{\lambda_m^{2i-2}} \left( \bar{I}_1^i - 3^i \right) + \frac{1}{D_1} \left( \frac{J^2 - 1}{2} - \ln J \right), \tag{3.7}$$

where $C_1 = \frac{1}{2}$, $C_2 = \frac{1}{20}$, $C_3 = \frac{11}{1050}$, $C_4 = \frac{19}{7000}$ and $C_5 = \frac{519}{673750}$ are constant coefficients, the symbols $\mu = 0.98765$, $\lambda_m = 7$ and $D_1 = 0.1$ are material parameters. Unless otherwise stated, all quantities in this work are dimensionless.

As shown in Fig. 3, uniaxial tension, uniaxial compression, equi-biaxial tension and equi-biaxial compression experiments are numerically performed on a $1 \times 1 \times 1$ unit cubic specimen, respectively. Based on the equations listed in Fig. 3, the corresponding values of $b_e$ and $\gamma$ of each type of experimental tests are stored in data sets $\mathcal{S}^{\mathrm{UT}} = \{(b_e^{\mathrm{UT}}, \gamma^{\mathrm{UT}})_{i=1,\ldots,5000}\}$, $\mathcal{S}^{\mathrm{UC}} = \{(b_e^{\mathrm{UC}}, \gamma^{\mathrm{UC}})_{i=1,\ldots,5000}\}$, $\mathcal{S}^{\mathrm{ET}} = \{(b_e^{\mathrm{ET}}, \gamma^{\mathrm{ET}})_{i=1,\ldots,5000}\}$, $\mathcal{S}^{\mathrm{EC}} = \{(b_e^{\mathrm{EC}}, \gamma^{\mathrm{EC}})_{i=1,\ldots,5000}\}$, respectively, where UT, UC, ET, and EC indicate uniaxial tension, uniaxial compression, equi-biaxial tension, and equi-biaxial compression tests, respectively. And the values of $J$ and $\sigma_m$ are the same in uniaxial and equi-biaxial tests, and are stored in data set $\mathcal{S} = \{(J, \sigma_m)_{i=1,\ldots,10000}\}$ without distinguishing the type of experimental test.

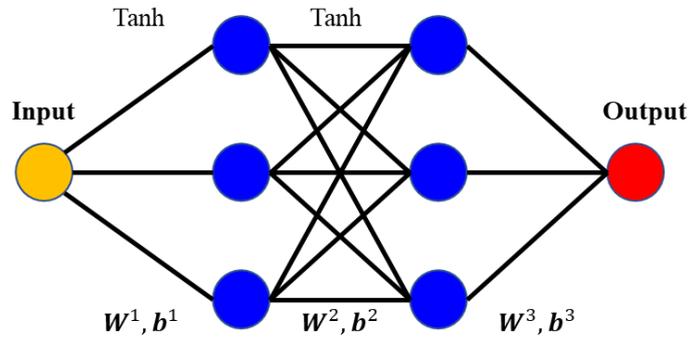

Fig. 4. The structure of the adopted artificial neural networks.

*3.1.3 An artificial neural network-based stress update strategy*

Since the data obtained by numerical experiments are discrete, smooth fitting processing is required to ensure the stability and robustness of the structural analysis and sensitivity analysis. Generally, fitting methods such as linear regression, polynomial regression, support vector regression, and artificial neural network (ANN) are often used to process discrete data. Thanks to the flexibility and strong generalization capabilities, five ANNs with the structure shown in Fig. 4 are adopted in this work. The inputs of the five ANNs are $b_e^{UT}$, $b_e^{UC}$, $b_e^{ET}$, $b_e^{EC}$ and $J$, and their corresponding outputs are $\gamma^{UT}$, $\gamma^{UC}$, $\gamma^{ET}$, $\gamma^{EC}$ and $\sigma_m$, respectively. In addition, in order to demonstrate the robustness of DDSA to data noise, 5% Gaussian random noise is added to each group of experimental data. The ideal data and the data with noise are shown in Fig. 5 and Fig. 6, which also illustrate that those ANNs achieve perfect smooth fitting even for the data with noise.

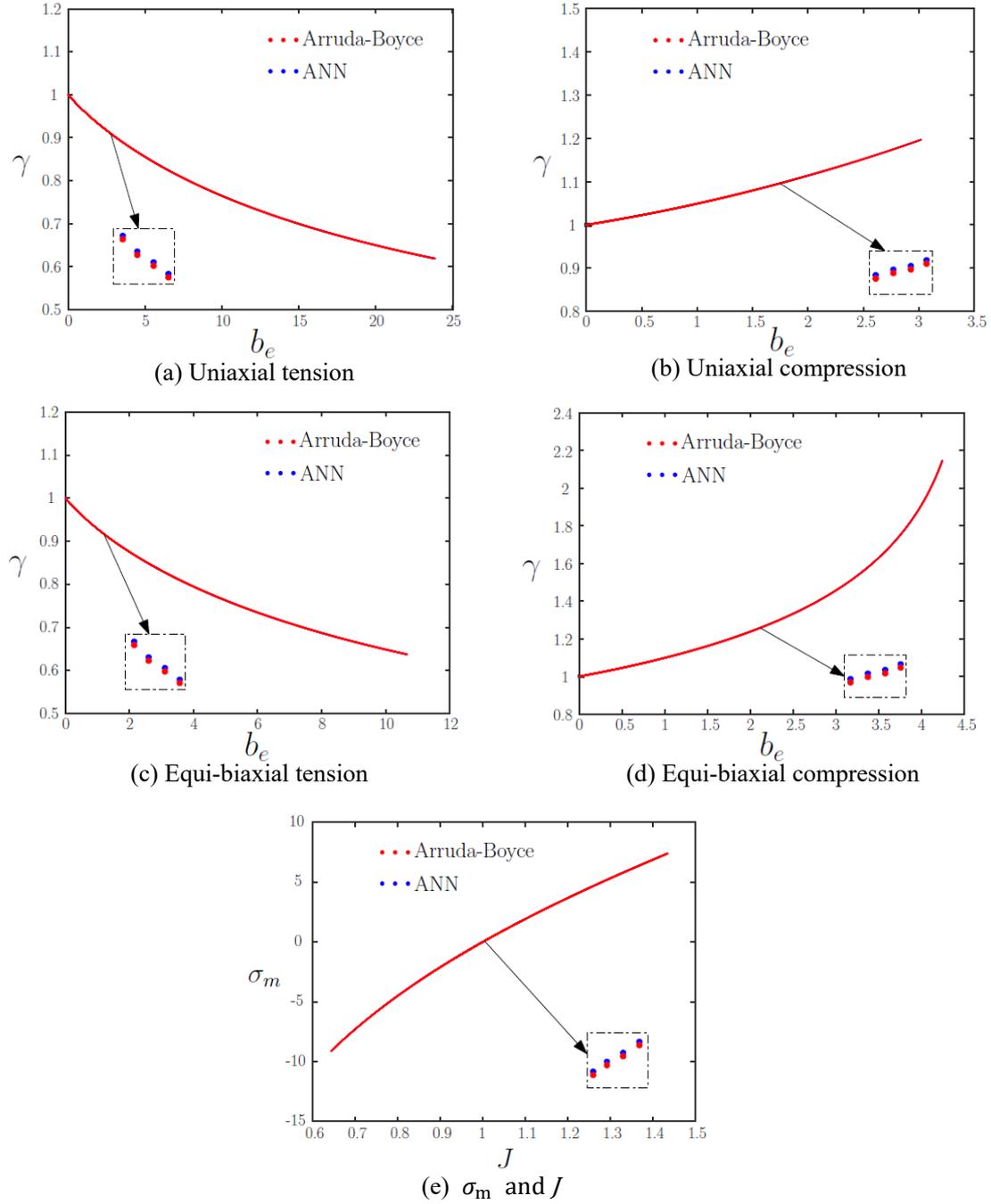

Fig. 5. Perfect data set of the Arruda-Boyce model fitted by the neural networks.

As we all know, in displacement-driven finite element analysis, $J$ and $b_e$ at any material point can be easily obtained by the displacement vector $\boldsymbol{U}$. Therefore, combining $J$ and the corresponding trained ANN, the spherical part $\sigma_m$ of the Cauchy stress can be uniquely determined. According to the trained ANNs, nevertheless, the $b_e$ of a material point corresponds to the four material parameters $\gamma^{UT}$, $\gamma^{UC}$, $\gamma^{ET}$, and $\gamma^{EC}$. Actually, the material parameter $\gamma$ is undoubtedly dependent on the stress state of the material point, which can be judged by the so-called the stress triaxiality introduced in [35]:

$$T = \frac{\sigma_m}{\sigma_e}. \tag{3.8}$$

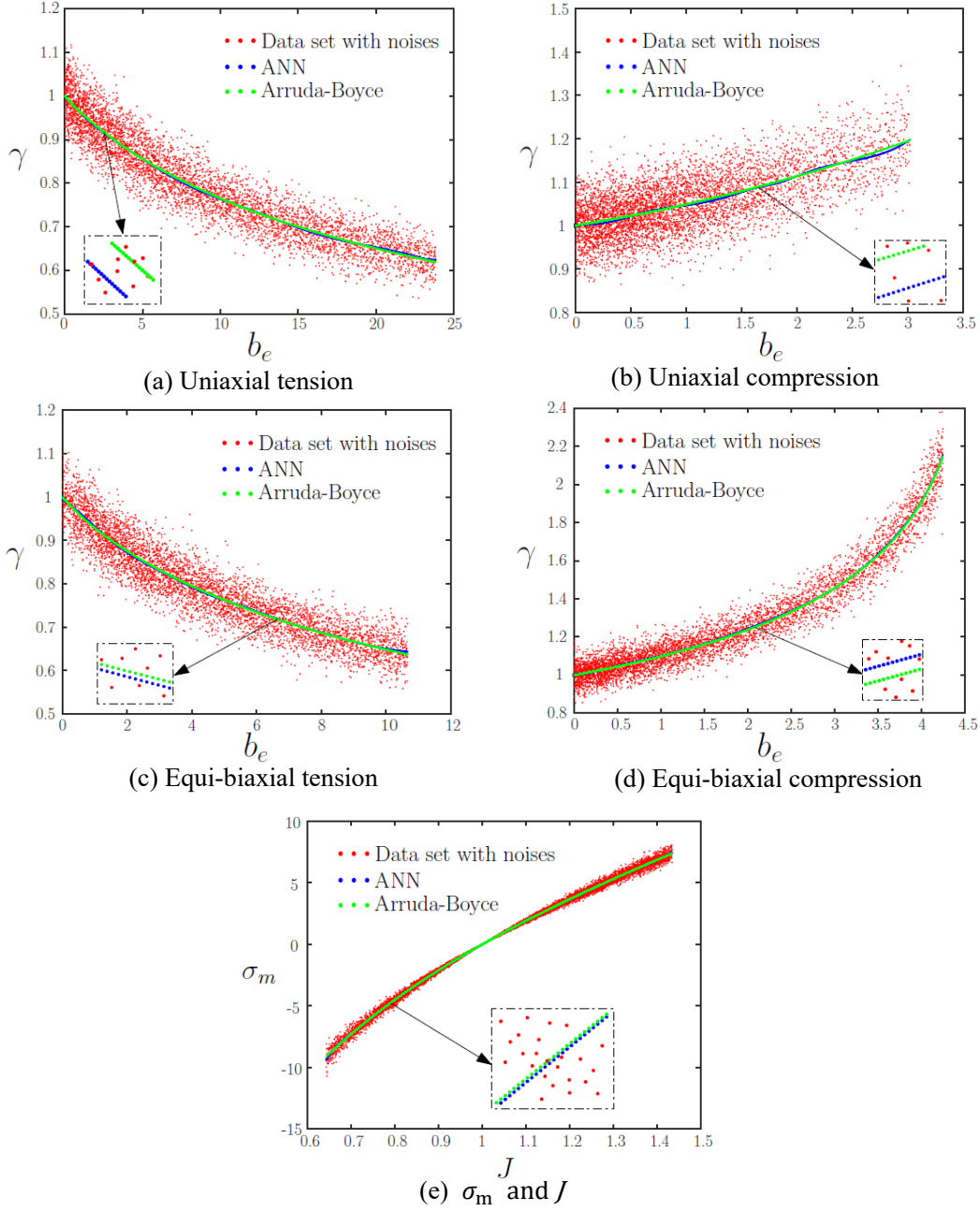

Fig. 6. Data set with noise of the Arruda-Boyce model fitted by the neural networks.

When the material point is in uniaxial tension and uniaxial compression states, the corresponding stress triaxiality is $\frac{1}{3}$ and $-\frac{1}{3}$, respectively; and when it is in equi-biaxial tension and equi-biaxial compression states, the corresponding stress triaxiality is $\frac{2}{3}$ and $-\frac{2}{3}$, respectively. Furthermore, we assume that the variation of $\sigma_e$ between the equivalent stress $\sigma_e^{UT}$ of uniaxial tension and the equivalent stress $\sigma_e^{ET}$ of equi-biaxial tension is linear. Then the value of $\sigma_e$ can be calculated as:

$$\sigma_e = (1-m)\sigma_e^{UT} + m\sigma_e^{ET}, \tag{3.9}$$

where $m = \frac{9T^2 - 3T\sqrt{12-27T^2}+2}{2(9T^2-1)}$. For more details of the derivation of $m$, please refer to the literature [35]. Substituting Eq. (3.9) into Eq. (3.6), the material coefficient $\gamma$ can be written as:

$$\gamma = (1-m)\gamma^{UT} + m\gamma^{ET}. \tag{3.10}$$

For the compression states of a material point, in the same way, the material coefficient can be obtained as:

$$\gamma = (1-m)\gamma^{UC} + m\gamma^{EC}. \tag{3.11}$$

Based on the above discussion, the determination of the material coefficient $\gamma$ can be summarized as:

$$\gamma = \begin{cases} \gamma^{EC}, & \text{if } T \leq -\frac{2}{3}, \\ (1-m)\gamma^{UC} + m\gamma^{EC}, & \text{if } -\frac{2}{3} < T < -\frac{1}{3}, \\ \gamma^{UC}, & \text{if } -\frac{1}{3} \leq T < 0, \\ \gamma^{UT}, & \text{if } 0 \leq T \leq \frac{1}{3}, \\ (1-m)\gamma^{UT} + m\gamma^{ET}, & \text{if } \frac{1}{3} < T < \frac{2}{3}, \\ \gamma^{ET} & \text{if } T \geq \frac{2}{3}. \end{cases} \tag{3.12}$$

Once the material coefficient $\gamma$ is given by Eq. (3.12), the Cauchy stress of the material point can be updated as:

$$\boldsymbol{\sigma} = \gamma \text{dev}(\boldsymbol{b}) + \sigma_m \boldsymbol{I}. \tag{3.13}$$

In the total Lagrangian formulation, the balance equation of the discrete form can be rewritten as:

$$\boldsymbol{R}(\boldsymbol{U}) = \boldsymbol{f}^{\text{int}} - \boldsymbol{f}^{\text{ext}} = \sum_{e=1}^{ne} \boldsymbol{L}_e^\top \boldsymbol{f}_e^{\text{int}} - \sum_{e=1}^{ne} \boldsymbol{L}_e^\top \boldsymbol{f}_e^{\text{ext}} = \boldsymbol{0}, \tag{3.14}$$

where the internal force vector $\boldsymbol{f}_e^{\text{int}} = \int_{^0\Omega_e} \boldsymbol{B}^\top (J\boldsymbol{\sigma} \cdot \boldsymbol{F}^{-\top}) \, \mathrm{d}^0 V_e$ and it can be calculated by combining Eq. (3.13). Then the analysis result can be obtained by the above-mentioned Newton-Raphson method. It is worth mentioning that the tangent stiffness matrix $\boldsymbol{K}_\mathrm{T}(\boldsymbol{U})$ used in the iteration process is updated by the finite difference method.

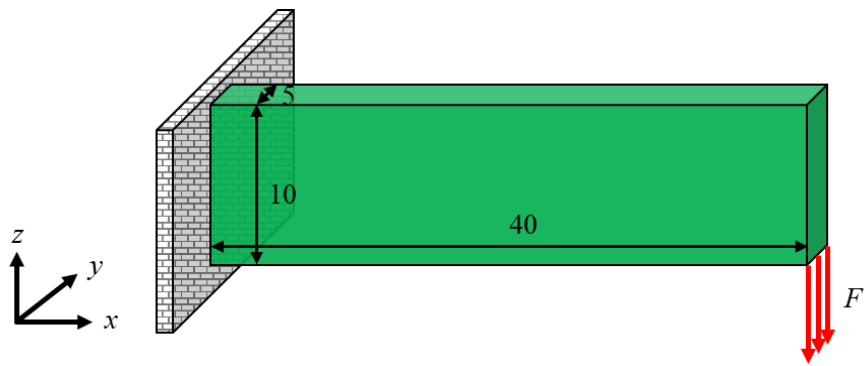

Fig. 7. The three-dimensional cantilever beam example

**3.2 Numerical verification of the improved DDSA algorithm**

As shown in Fig. 7, a cantilever beam structure with a size of $40 \times 5 \times 10$ is fixed on the left and discretized by $80 \times 10 \times 20$ eight-node hexahedral elements. Uniform vertical downward loads are applied at all the nodes located at the lower edge of the right face for three cases with increasing amplitudes, i.e., $F = 0.03$, $F = 0.05$, $F = 0.07$, respectively. Based on the reference Arruda-Boyce model mentioned above and the constructed database, the analysis results of the three cases are shown in Fig. 8, which includes the Z-direction displacement contour, the maximum displacement of Z-direction and the structural compliance. Furthermore, for a clearer comparison, the relative errors of the results of the DDSA algorithm and their references are listed in Table 1 (all smaller than 2.5%). It can be found that no matter whether noise exists or nor, the DDSA algorithm can accurately obtain the analysis results with very small errors. This example illustrates the effectiveness of the proposed DDSA algorithm and its high robustness to noise.

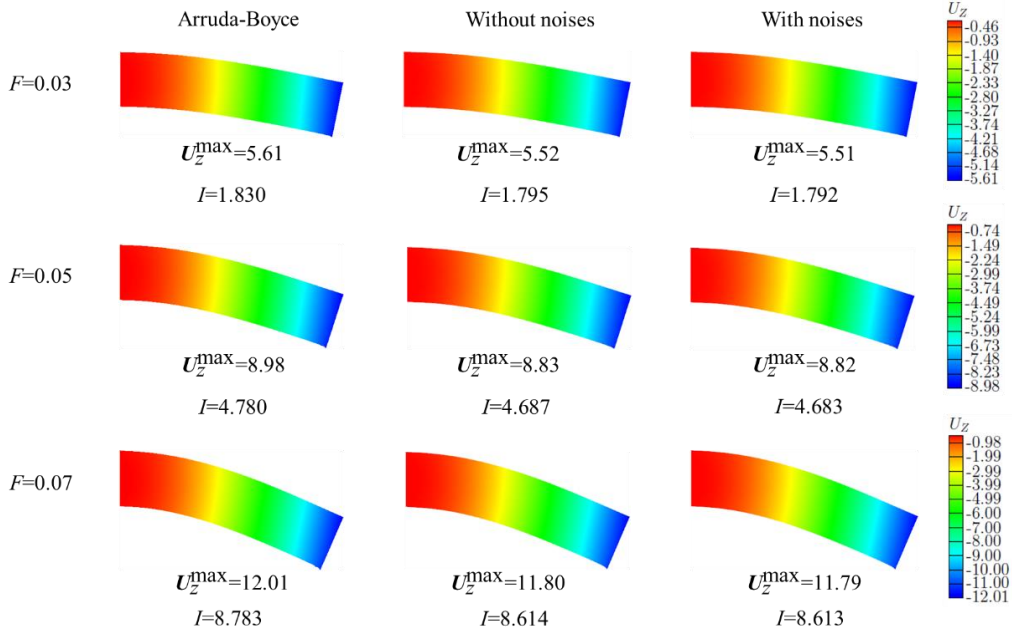

Fig. 8. Comparison of the analysis results obtained by the constitutive model-based analysis and improved DDSA algorithms.

Table 1 The relative errors of the results of the DDSA and their reference values.

| Force | $U_{\mathrm{Z}}^{\max}$ (Without noise) | $U_{\mathrm{Z}}^{\max}$ (With noise) | Compliance (Without noise) | Compliance (With noise) |
| --- | --- | --- | --- | --- |
| 0.03 | 1.60% | 1.78% | 1.91% | 2.08% |
| 0.05 | 1.67% | 1.78% | 1.95% | 2.03% |
| 0.07 | 1.75% | 1.83% | 1.92% | 1.94% |

**3.3 Sensitivity analysis and material interpolation model of the DDTO framework**

Since the DDSA method is under the displacement-driven finite element analysis framework, the problem formulation of DDTO for three-dimensional hyperelastic continuum structures under finite deformation is the same in form as Eq. (2.3), which can be formulated as

Find $\quad \boldsymbol{D} = \left((\boldsymbol{D}^1)^\top, \ldots, (\boldsymbol{D}^i)^\top, \ldots, (\boldsymbol{D}^{nv})^\top\right)^\top, {}^0\boldsymbol{u}(\boldsymbol{D}; {}^0\boldsymbol{x})$,

Minimize $\quad I(\boldsymbol{D}, {}^0\boldsymbol{u}) = \int_{{}^0\Omega} H_\epsilon^\alpha\left(\chi^s(\boldsymbol{D}; {}^0\boldsymbol{x})\right) {}^0\boldsymbol{f} \cdot {}^0\boldsymbol{u}\, \mathrm{d}^0\mathrm{V} + \int_{{}^0S_t} {}^0\boldsymbol{t} \cdot {}^0\boldsymbol{u}\, \mathrm{d}^0\mathrm{S}$,

S. t.

$$\int_{^0\Omega} H_\epsilon^\alpha(\chi^s(\boldsymbol{D}; {}^0\boldsymbol{x}))\delta\boldsymbol{F}^\top : \boldsymbol{P} \mathrm{d}^0\mathrm{V} = \int_{^0\Omega} H_\epsilon^\alpha(\chi^s(\boldsymbol{D}; {}^0\boldsymbol{x})){}^0\boldsymbol{f} \cdot \delta^0\boldsymbol{u} \mathrm{d}^0\mathrm{V} + \int_{^0S_t} {}^0\boldsymbol{t} \cdot \delta^0\boldsymbol{u} \mathrm{d}^0\mathrm{S}, \forall \delta^0\boldsymbol{u} \in \mathcal{U}_{\mathrm{ad}}^0,$$

$$\boldsymbol{P} = \mathcal{F}(\boldsymbol{F}, J; \boldsymbol{NN}_1, \boldsymbol{NN}_2, \boldsymbol{NN}_3, \boldsymbol{NN}_4, \boldsymbol{NN}_5),$$

$$\int_{^0\Omega} H_\epsilon^\alpha(\chi^s(\boldsymbol{D}; {}^0\boldsymbol{x})) \mathrm{d}^0\mathrm{V} \leq {}^0\bar{V},$$

$$\boldsymbol{D} \subset \mathcal{U}_D,$$

$${}^0\boldsymbol{u} = \boldsymbol{0}, \mathrm{on}\, {}^0S_u, \tag{3.15}$$

where $\boldsymbol{NN}_1$, $\boldsymbol{NN}_2$, $\boldsymbol{NN}_3$, $\boldsymbol{NN}_4$, and $\boldsymbol{NN}_5$ are the network parameters corresponding to the five ANNs trained in sub-section 3.1.3, respectively. *It is worth noting that for any material point in the structure, the 1st Piola–Kirchhoff stress tensor $\boldsymbol{P}$ is no longer determined by the analytic-function-based $W$, but is given by combining the constructed database and the stress update strategy in the DDTO framework.*

The solution strategies of DDTO framework, such as sensitivity analysis and numerical stabilization technique, can perfectly transplant their counterparts CMTO uses, e.g., the adjoint sensitivity analysis method and the redundant degrees of freedom removal technique in subsection 2.3. As a result, the sensitivity of the objective function in DDTO framework is the same in form as Eq. (2.11). As compared to the sensitivity analysis of CMTO, the only difference in sensitivity of the objective function in DDTO method is that the nodal internal force vector $\tilde{\boldsymbol{f}}_e^{\mathrm{int}}$ of the $e$th element is calculated by the stress update strategy and the tangent stiffness matrix $\boldsymbol{K}_\mathrm{T}^E(\boldsymbol{D}, \widetilde{\boldsymbol{U}}^E)$ used in Eq. (2.12) is given by the finite difference method.

Additionally, according to Eqs. (2.5), (3.2), (3.3a) and (3.3b), the material interpolation model of DDTO can be expressed as:

$$\boldsymbol{\sigma}_e = \frac{\sum_{i=1}^{nn} H_\epsilon^\alpha(\chi^s(\boldsymbol{D}; {}^0\boldsymbol{x})_i^e)}{nn}(\gamma\mathrm{dev}(\boldsymbol{b}) + \sigma_\mathrm{m}\boldsymbol{I}) = \frac{\sum_{i=1}^{nn} H_\epsilon^\alpha(\chi^s(\boldsymbol{D}; {}^0\boldsymbol{x})_i^e)}{nn}\boldsymbol{\sigma}_0, \tag{3.16}$$

where $\boldsymbol{\sigma}_e$ is the Cauchy stress of the $e$th element, and $\boldsymbol{\sigma}_0$ is the Cauchy stress obtained by combining the constructed database and the stress update strategy (in sub-section 3.1.2 and 3.1.3). The flowchart of the proposed DDTO framework is shown in Fig. 9.

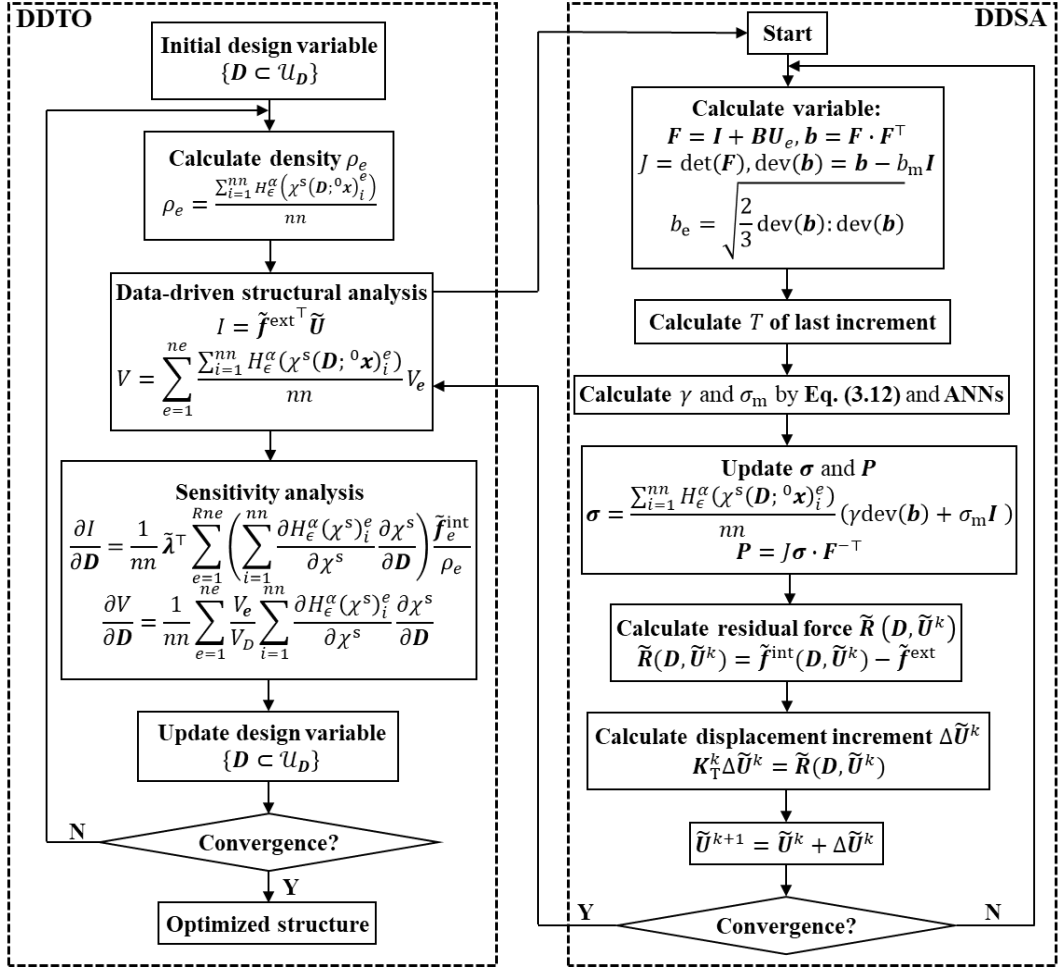

Fig. 9. The flowchart of the DDTO framework.

## 4. Numerical examples

In this section, several examples are studied to demonstrate the effectiveness and stability of the proposed DDTO framework for three-dimensional continuum structures under finite deformation. The experimental data and the analytic-function based constitutive model used in the numerical examples are given in sub-section 3.1.2 unless otherwise stated. The optimization problems are solved by the well-known Method of Moving Asymptotes (MMA) [39] and terminated by the convergence criterion of $V_{\text{var}} < 10^{-3}$ and $obj_{\text{var}} < 10^{-3}$, where $V_{\text{var}}$, $obj_{\text{var}}$ are defined as

$$V_{\text{var}} = \left| \frac{V(\boldsymbol{D})}{V_D} - \frac{^0\bar{V}}{V_D} \right|, \quad obj_{\text{var}} = \begin{cases} 1.0, & \text{if } iter < 5, \\ \frac{\sum_{i=iter-4}^{iter} |obj^i - obj^{iter}|}{5 obj^{iter}}, & \text{otherwise,} \end{cases} \quad (4.1)$$

where $V(\boldsymbol{D})$ is the volume of the current design, $iter$ is the number of current iteration, $obj^i$ is the objective function value at the $i$th iteration.

**4.1 The three-dimensional cantilever beam example**

In order to verify the effectiveness of the 3D-MMV based explicit topology optimization algorithm for continuum structures with an analytic-function based constitutive model (i.e., the CMTO method), the cantilever beam example described by Section 3.2 and Fig. 7 is investigated first. The amplitudes of the external nodal loads are set as $F = 0.02$, $F = 0.025$, $F = 0.04$, $F = 0.06$, respectively. Setting the initial design shown in Fig. 10 and the upper bound of the volume fraction of the solid material as 0.5, the optimized results of the four loading cases are shown in Table 2, respectively. According to the optimized results in Table 2, significant differences are found between the optimized structures of the four cases, which are consistent with the results in literature [8]. Furthermore, the stable iteration history curve of the case with $F = 0.06$ is presented in Fig. 11.

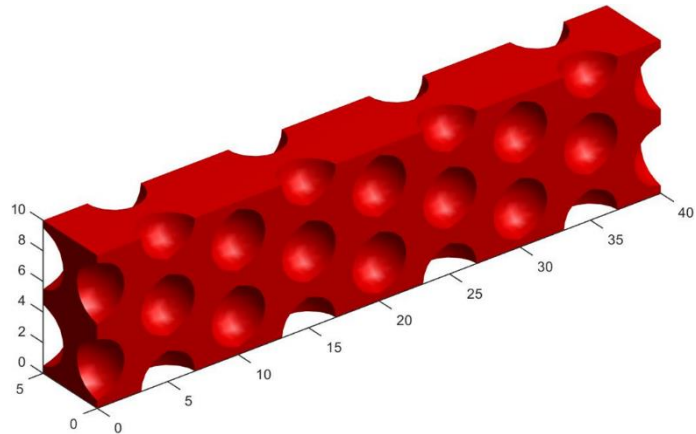

Fig. 10. The initial design for the three-dimensional cantilever beam.

Table 2 The optimized results of the four cases solved by the CMTO method.

| Force | Optimized structure | Compliance |
|---|---|---|
| 0.02 | 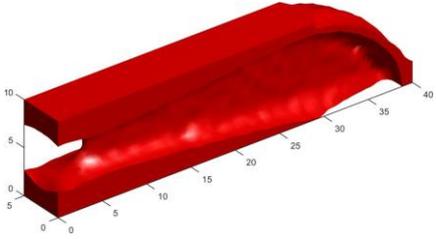 | 1.06 |
| 0.025 | 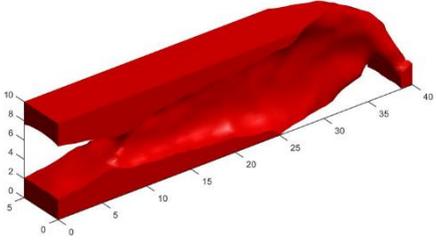 | 1.64 |
| 0.04 | 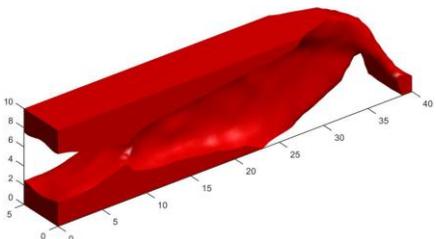 | 3.86 |
| 0.06 | 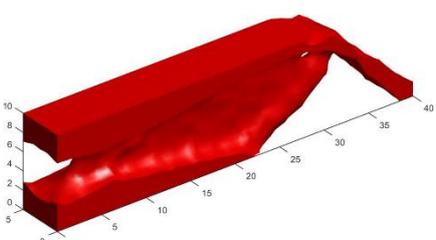 | 7.74 |

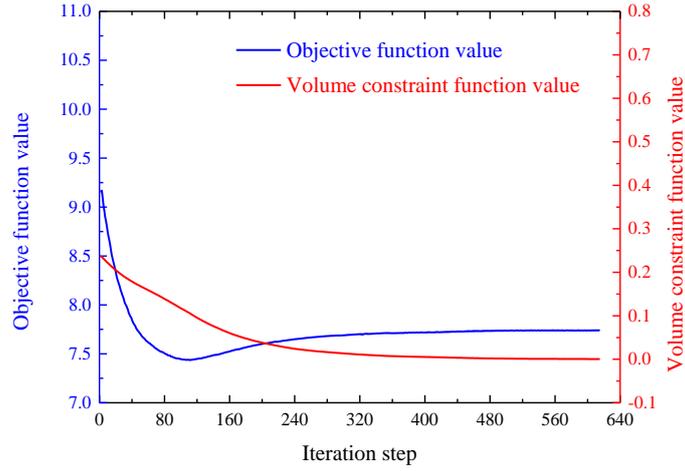

Fig. 11. The iterative history curve of $F = 0.06$ of the cantilever beam example.

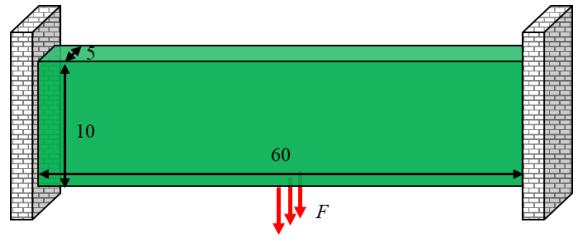

Fig. 12. The three-dimensional two-ends clamped beam example.

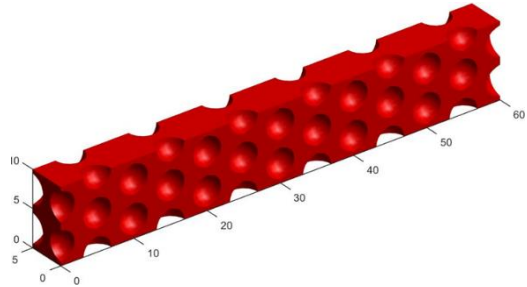

Fig. 13. The initial design for the three-dimensional two-ends clamped beam.

**4.2 The three-dimensional two-ends clamped beam example**

As shown in Fig. 12, the design domain with a size of $60 \times 5 \times 10$ is discretized by $120 \times 10 \times 20$ eight-node hexahedra elements. Both ends of the design domain are fixed and the downward nodal loads of three cases, i.e., $F = 0.1$, $F = 0.3$, $F = 0.6$, respectively, are imposed on all the nodes located at the center line of the bottom of the design domain. The initial design shown in Fig. 13 contains 64 voids with 4160 design variables and the upper bound of the volume fraction of the solid material is set as 0.5. This example is solved both by the proposed CMTO method and the DDTO method with the dataset shown in Fig. 5.

The optimized results corresponding to the three cases and the iterative curves of $F = 0.1$ are shown in Table 3 and Fig. 14, respectively. It can be found that, both of the proposed algorithms converge stably. Furthermore, as illustrated by Table 3, both the optimized structures of DDTO method and CMTO method have significant topology changes as the variation of the amplitude of the external load, which is caused by the consideration of the material and geometrical nonlinearities. Besides, only small differences are observed between the optimized structures obtained by the DDTO and CMTO methods, and the relative errors of the objective function values for the three load cases are all smaller than 3%. These results demonstrate the correctness and validity of the proposed DDTO method.

Table 3 The optimized results of the three cases solved by the CMTO and DDTO methods.

| Framework | Force | Optimized structure | Compliance | Relative error |
|---|---|---|---|---|
| CMTO | 0.1 | 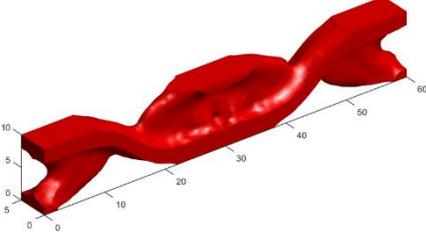 | 2.39 | |
| DDTO | 0.1 | 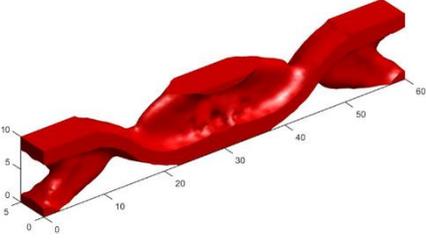 | 2.32 | 2.93% |
| CMTO | 0.3 | 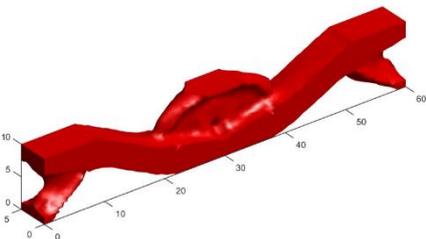 | 18.42 | |
| DDTO | 0.3 | 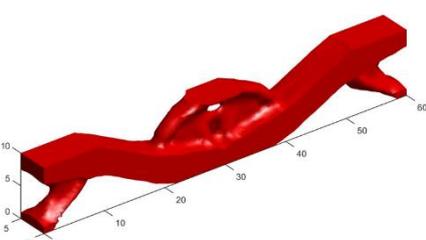 | 17.89 | 2.88% |

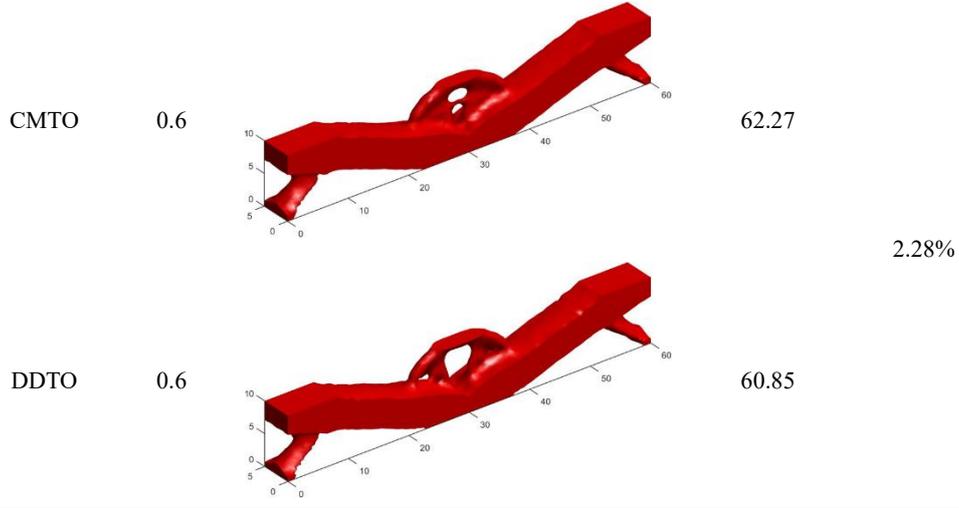

| | | | |
|---|---|---|---|
| CMTO | 0.6 | | 62.27 |
| | | | 2.28% |
| DDTO | 0.6 | | 60.85 |

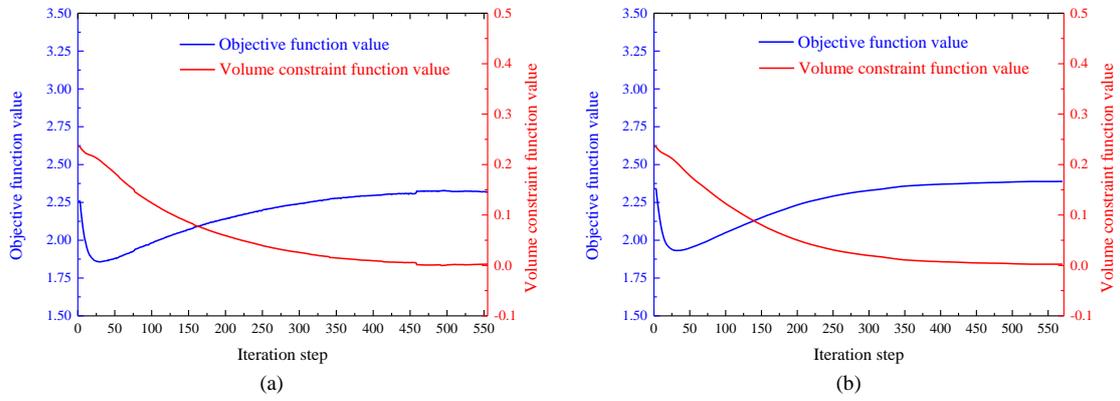

Fig. 14. The iterative history curves of $F = 0.1$, (a) the DDTO method, (b) the CMTO method.

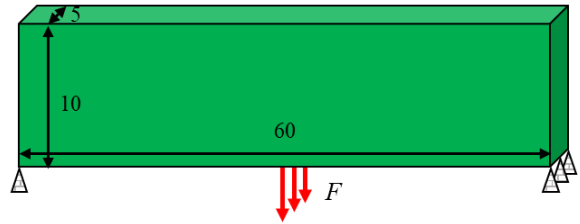

Fig. 15. The three-dimensional simply-supported beam example.

## 4.3 The three-dimensional simply-supported beam example

This simply-supported beam example shown in Fig. 15 is to test the robustness of the proposed DDTO algorithm against data noise. The size of the design domain is $60 \times 5 \times 10$ and $120 \times 10 \times 20$ uniform eight-node hexahedra elements are used for discretizing the design domain. The vertical downward nodal loads with magnitudes of $F = 0.1$ are imposed on all the nodes located at the center line of the bottom of the design domain. The initial design is the same as in example 4.2, as shown in Fig. 13. The upper bound of the volume fraction of the base material is set as 0.5. Under the conditions mentioned above, this example is solved by the DDTO algorithm with noisy data set (as shown in Fig. 6), the DDTO algorithm with ideal data set and the CMTO algorithm, respectively. The corresponding optimized structures, compliance values and their relative errors are shown in Table 4. Those optimized designs are very close and the relative error of the compliance values of the DDTO designs are smaller than 3% as compared to the results obtained by the CMTO method, which validates that the proposed DDTO framework is also robust to data noise.

Table 4 The optimized results of the three-dimensional simply-supported beam.

| Framework | Optimized structure | Compliance | Relative error |
|---|---|---|---|
| CMTO | 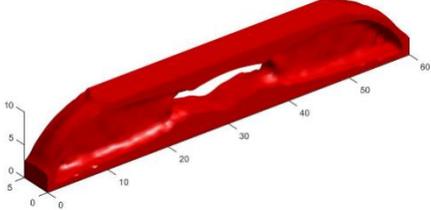 | 7.28 | —— |
| DDTO (Without noise) | 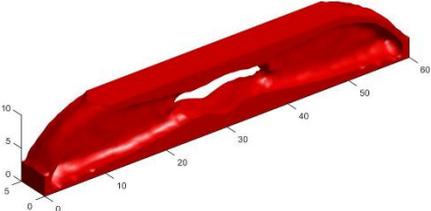 | 7.09 | 2.61% |
| DDTO (With noise) | 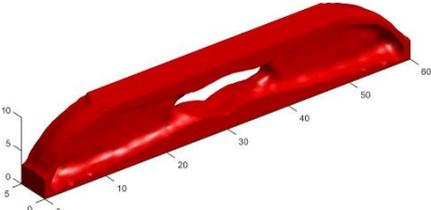 | 7.08 | 2.75% |

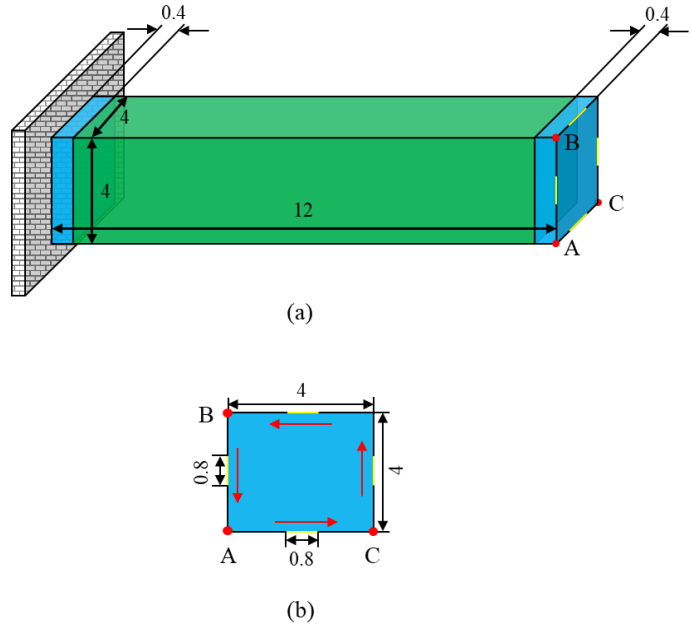

(a)

(b)

Fig. 16. The three-dimensional torsion beam example, (a) the design domain with left-hand side fixed, (b) uniform loads are applied to the four yellow lines on the face BAC, resulting in a counterclockwise torque.

**4.4 The three-dimensional torsion beam example**

Finally, the three-dimensional torsion beam is discussed to verify the effectiveness of the proposed DDTO algorithm for optimized structures with complex stress states. As illustrated in Fig. 16(a), a $12 \times 4 \times 4$ beam is fixed at the left-hand side and the nodal loads of magnitude 0.02 are applied at all the nodes located at the four yellow lines on the right face BAC, resulting in a counter-clockwise torque (see Fig. 16(b) for reference). The left and right sides (both have a dimension of $0.4 \times 4 \times 4$) are defined as non-designable solid domain. For the DDSA, $60 \times 20 \times 20$ eight-node hexahedra elements are created. The initial design shown in Fig.17 contains 88 voids with 5720 design variables, and the upper bounds of the volume fraction of the material are set as $^0\bar{V} = 0.5$ and $^0\bar{V} = 0.3$, respectively. This example is solved by the proposed DDTO method and the CMTO algorithm as well.

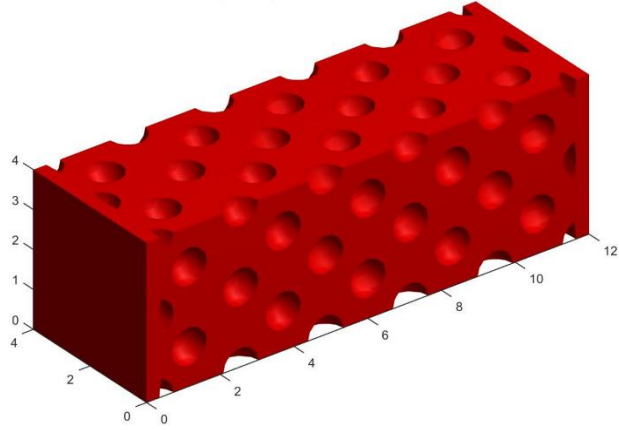

Fig. 17. The initial design for the three-dimensional torsion beam.

The optimized results and the iterative history curves of $^0\bar{V} = 0.5$ are shown in Table 5, Table 6 and Fig. 18, respectively. As shown in Tables 5 and 6, both algorithms obtain a thin-walled cylindrical structure, and the cylindrical wall becomes thinner as the upper bound of the volume fraction of the base material becomes smaller. Furthermore, the relative errors of the compliance values are both small for the cases $^0\bar{V} = 0.5$ and $^0\bar{V} = 0.3$, i.e., 2.94% and 2.17% respectively, which shows that the proposed DDTO method is also effective despite complex stress states.

Thanks to the adoption of the redundant degrees of freedom removal technique, the DDTO algorithm can not only achieve a stable convergence during the optimization process (see Fig. 18), but also significantly reduce the number of involved DOFs in the nonlinear finite element analysis as shown in Fig. 19. In addition, the optimized structure can be directly constructed by the CAD system as shown in Fig. 20, since the topology of the structure is described by explicit geometric parameters under the proposed DDTO framework.

Table 5 The optimized results of the three-dimensional torsion beam of $^0\bar{V} = 0.5$.

| | CMTO | DDTO |
|---|---|---|
| Optimized structure | | |
| Cross section | 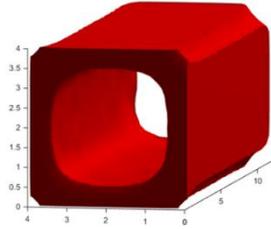 | 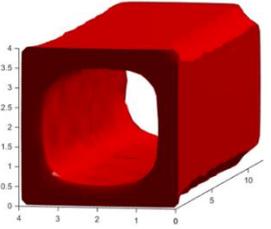 |
| Compliance | 0.34 | 0.33 |
| Relative error | ―― | 2.94% |

Table 6 The optimized results of the three-dimensional torsion beam of $^0\bar{V} = 0.3$.

| | CMTO | DDTO |
|---|---|---|
| Optimized structure | 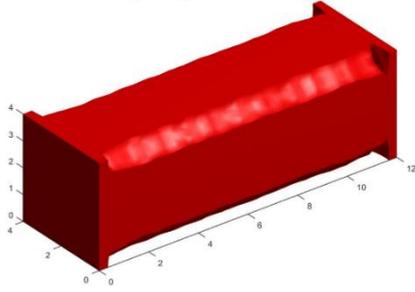 | 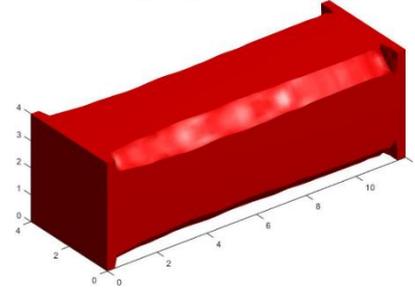 |
| Cross section | 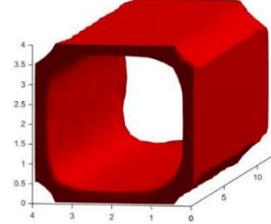 | 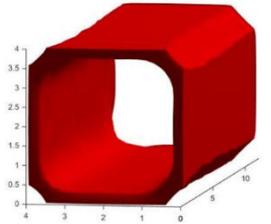 |
| Compliance | 0.46 | 0.45 |
| Relative error | ―― | 2.17% |

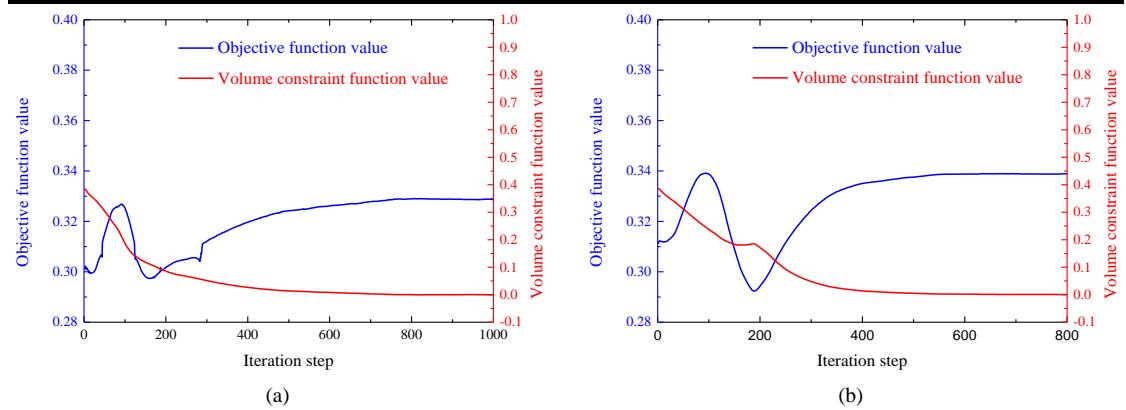

(a)          (b)

Fig. 18. The iterative history curves of the (a) DDTO algorithm and (b) CMTO algorithm.

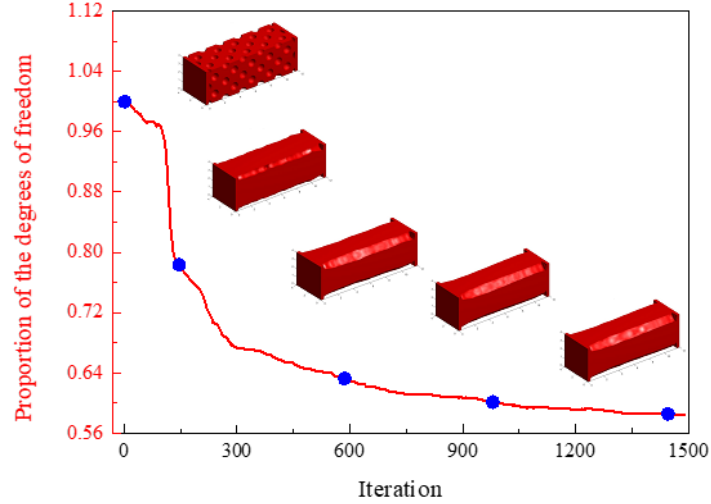

Fig. 19. Iteration history of the proportion of retained degrees of freedom for structural analysis of DDTO with $^0\bar{V} = 0.3$.

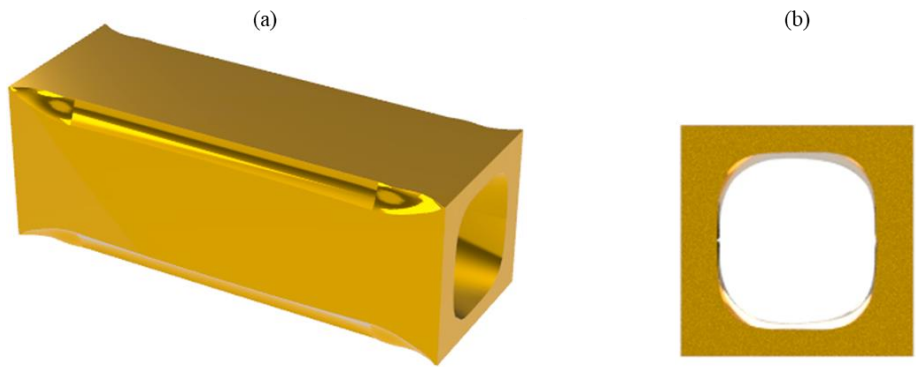

Fig. 20. CAD model of the optimized torsion beam of $^0\bar{V} = 0.5$ without non-designable solid domain, (a) global perspective, (b) sectional perspective.

## 5. Concluding remarks

In the present work, a data-driven topology optimization (DDTO) framework for three-dimensional continuum structures under finite deformation is proposed. The advantages of proposed DDTO framework can be summarized as follows: (1) In the DDTO framework, topology optimization of the three-dimensional continuum structure under finite deformation is implemented only by the uniaxial and equi-biaxial experimental data, without using the analytic-function based constitutive models. (2) Since the data-driven analysis method is still under the framework of the finite element method, the solution strategies used in traditional topology optimization method such as sensitivity analysis and redundant degrees of freedom removal technique can be perfectly transplanted to the DDTO framework. (3) Described by a series of NURBS surfaces, the optimized structure has explicit geometry and can be constructed in CAD system directly. Although only minimum end compliance design is considered in the present work, with the assumption of hyperelasticity of the considered material, the DDTO framework can also be extended to topology optimization problems considering the effects of multi-material, stress, frequency, buckling, etc. For materials with more general nonlinear behaviors, e.g., elastoplastic and viscoelastic materials, the data-driven topology optimization of the three-dimensional continuum structure still requires further efforts.

**Acknowledgment**

The financial supports from the National Natural Science Foundation (11821202, 11732004, 12002073, 11872139), the National Key Research and Development Plan (2020YFB1709401), Dalian Talent Innovation Program (2020RQ099), and 111 Project (B14013) are gratefully acknowledged.

**Conflict of interest**

The authors declare they have no conflict of interest.

**Replication of results**

Code and data for replication can be provided up on request.

**Appendix A**

$$r(\theta, \varphi) = \|S(u(\theta), v(\varphi))\| = \left\|\sum_{k=1}^{K}\sum_{l=1}^{L} N_{k,p}(u(\theta)) N_{l,q}(v(\varphi)) \boldsymbol{P}^{k,l}\right\|, \tag{A.1}$$

where $N_{k,p}(u), N_{l,q}(v)$ are the B-spline basis functions of $p$th and $q$th orders, $K$ and $L$ are the total number of control points in the $\theta$ and $\varphi$ directions, respectively. The symbols $\boldsymbol{P}^{k,l} = \left(P_x^{k,l}, P_y^{k,l}, P_z^{k,l}\right)^\mathsf{T}, k = 0, \dots, K, l = 0, \dots, L$ are the coordinates of the corresponding control points. In order to avoid self-intersection of the surfaces, the control points are constructed as follows:

$$P_x^{k,l} = r^{k,l} \sin(\varphi_l) \cos(\theta_k), \tag{A.2a}$$

$$P_y^{k,l} = r^{k,l} \sin(\varphi_l) \sin(\theta_k), \tag{A.2b}$$

$$P_z^{k,l} = r^{k,l} \cos(\varphi_l), \tag{A.2c}$$

where $r^{k,l}, k = 0, \dots, K, l = 0, \dots, L$ is the distance from the control point to the center of the void, $\theta_k = \frac{2k\pi}{K}, k = 0, \dots, K, \varphi_l = \frac{l\pi}{L}, l = 0, \dots, L$. Since the surface is a closed region, it is also required that

$$r^{0,l} = r^{K,l}, \quad l = 0, \dots, L, \tag{A.3a}$$

$$r^{k,0} = r^{0,0}, \quad k = 1, \dots, K, \tag{A.3b}$$

$$r^{k,L} = r^{0,L}, \quad k = 1, \dots, K. \tag{A.3c}$$

**Appendix B**

$$\frac{\partial \chi^s}{\partial \boldsymbol{D}} = \sum_{i=1}^{nv} \frac{\partial \chi^s}{\partial \chi^i(\boldsymbol{D}^i; \boldsymbol{x})} \frac{\partial \chi^i(\boldsymbol{D}^i; \boldsymbol{x})}{\partial \boldsymbol{D}} = \sum_{i=1}^{nv} \frac{e^{-\lambda \chi^i(\boldsymbol{D}^i; \boldsymbol{x})}}{\sum_{i=1}^{nv} e^{-\lambda \chi^i(\boldsymbol{D}^i; \boldsymbol{x})}} \frac{\partial \chi^i(\boldsymbol{D}^i; \boldsymbol{x})}{\partial \boldsymbol{D}}, \tag{B.1}$$

In Eq. (2.11), the derivative of $\chi^i(\boldsymbol{D}^i; \boldsymbol{x})$ with respect to each design variable can be calculated as follows:

$$\frac{\partial \chi^i(\boldsymbol{D}^i; \boldsymbol{x})}{\partial C_x^i} = \frac{C_x^i}{\sqrt{(x - C_x^i)^2 + (y - C_y^i)^2 + (z - C_z^i)^2}}, \tag{B.2}$$

$$\frac{\partial \chi^i(\boldsymbol{D}^i; \boldsymbol{x})}{\partial C_y^i} = \frac{C_y^i}{\sqrt{(x - C_x^i)^2 + (y - C_y^i)^2 + (z - C_z^i)^2}}, \tag{B.3}$$

$$\frac{\partial \chi^i(\boldsymbol{D}^i; \boldsymbol{x})}{\partial C_z^i} = \frac{C_z^i}{\sqrt{(x - C_x^i)^2 + (y - C_y^i)^2 + (z - C_z^i)^2}}, \tag{B.4}$$

$$\frac{\partial \chi^i(\boldsymbol{D}^i; \boldsymbol{x})}{\partial r^{k,l}} = \frac{1}{r(\theta, \varphi)} \left[ S_x(\theta, \varphi) \frac{\partial S_x(\theta, \varphi)}{\partial r^{k,l}} + S_y(\theta, \varphi) \frac{\partial S_y(\theta, \varphi)}{\partial r^{k,l}} + S_z(\theta, \varphi) \frac{\partial S_z(\theta, \varphi)}{\partial r^{k,l}} \right],$$

$$k = 0, \dots, K, l = 0, \dots, L, \tag{B.5}$$

where

$$\frac{\partial S_x(\theta, \varphi)}{\partial r^{k,l}} = \sum_{n=1}^{K} \sum_{m=1}^{L} N_{n,p}(u(\theta)) N_{m,q}(v(\varphi)) \frac{\partial P_x^{n,m}}{\partial r^{k,l}}, n = 0, \dots, K, m = 0, \dots, L, \tag{B.6}$$

$$\frac{\partial S_y(\theta, \varphi)}{\partial r^{k,l}} = \sum_{n=1}^{K} \sum_{m=1}^{L} N_{n,p}(u(\theta)) N_{m,q}(v(\varphi)) \frac{\partial P_y^{n,m}}{\partial r^{k,l}}, n = 0, \dots, K, m = 0, \dots, L, \tag{B.7}$$

$$\frac{\partial S_z(\theta, \varphi)}{\partial r^{k,l}} = \sum_{n=1}^{K} \sum_{m=1}^{L} N_{n,p}(u(\theta)) N_{m,q}(v(\varphi)) \frac{\partial P_z^{n,m}}{\partial r^{k,l}}, n = 0, \dots, K, m = 0, \dots, L. \tag{B.8}$$

If $n = k$ and $m = l$, $\partial P_x^{k,l}/\partial r^{k,l}$, $\partial P_y^{k,l}/\partial r^{k,l}$, $\partial P_z^{k,l}/\partial r^{k,l}$ are be expressed as

$$\frac{\partial P_x^{k,l}}{\partial r^{k,l}} = \sin(\varphi_l)\cos(\theta_k), \quad \frac{\partial P_y^{k,l}}{\partial r^{k,l}} = \sin(\varphi_k)\sin(\theta_l), \quad \frac{\partial P_z^{k,l}}{\partial r^{k,l}} = \cos(\varphi_l), \tag{B.9}$$

otherwise, $\partial P_x^{n,m}/\partial r^{k,l} = 0$, $\partial P_y^{n,m}/\partial r^{k,l} = 0$, $\partial P_z^{n,m}/\partial r^{k,l} = 0$.